\documentclass{amsart}
\usepackage{setspace}

\theoremstyle{plain}
\newtheorem{theorem}{Theorem}[section]
\newtheorem{proposition}{Proposition}[section]
\newtheorem{conjecture}{Conjecture}[section]

\theoremstyle{definition}

\newtheorem{cor}{Corollary}
\newtheorem{example}{Example}[section]
\newtheorem{exercise}{Exercise}

\newenvironment{renumerate}%
{%
\begin{enumerate}}%
{\end{enumerate}%
}%

\newenvironment{theo}[1]%
{\begin{theorem}\label{#1}}%
{\end{theorem}}

{\begin{proposition}\label{#1}}%
{\end{proposition}}

{\begin{conjecture}\label{#1}}%
{\end{conjecture}}

\newenvironment{definition}%
{\vskip6pt%
\noindent%
{\bf Definition.}}%
{\vskip6pt}

{\vskip6pt%
\noindent%
{\it Remark.}}%
{\vskip6pt}

{\vskip6pt%
\noindent%
{\it Remarks}. %
\begin{renumerate}}%
{\end{renumerate}\vskip6pt}

{\begin{exercise}\label{#1}}%
{\end{exercise}}

\newenvironment{ex}[1]%
{\begin{example}\label{#1}}%
{\end{example}}

\newcommand{\R}{\text{${\mathbb R}$}}
\newcommand{\C}{\text{$\mathbb C$}}

\newcommand{\Z}{\text{$\mathbb Z$}}

\newcommand{\J}{\text{$\mathcal{J}$}}

\newcommand{\gO}{\text{$\Omega$}}
\newcommand{\go}{\text{$\omega$}}

\newcommand{\del}{\text{$\partial$}}

\newcommand{\tensor}{\otimes}

\newcommand{\mc}[1]{\text{$\mathcal{#1}$}}

\newcommand{\into}{\longrightarrow}
\newcommand{\noqed}{\let\qed\relax}

\newcommand{\IP}[1]{\langle #1 \rangle}
\newcommand{\gcx}{generalized complex}
\newcommand{\gcs}{generalized complex structure}

\newcommand{\gcm}{generalized complex manifold}
\newcommand{\gcms}{generalized complex manifolds}

\numberwithin{equation}{section}

\author{Gil R. Cavalcanti} 
\author{Marco Gualtieri}
\email{gilrc@maths.ox.ac.uk}
\email{mgualt@mit.edu}

\title{A surgery for generalized complex structures on 4-Manifolds}
\begin{document}
\maketitle
\begin{abstract}
We introduce a surgery for \gcms\ whose input is a symplectic
4-manifold containing a symplectic 2-torus with trivial normal
bundle and whose output is a 4-manifold endowed with a \gcs\
exhibiting type change along a 2-torus.  Performing this surgery on
a $K3$ surface, we obtain a \gcs\ on $3 \C P^2 \# 19 \overline{\C
P^2}$, which has vanishing Seiberg--Witten invariants and hence does
not admit complex or symplectic structure.
\end{abstract}

\markboth{\sc G. R. Cavalcanti and M. Gualtieri}{\sc A surgery for
generalized complex 4-manifolds}

\section*{Introduction}
Generalized complex structures, introduced by Hitchin \cite{Hi03}
and developed by the second author in \cite{Gu03}, are a
simultaneous generalization of complex and symplectic structures. In
this paper we answer, in the affirmative, the question of whether
there exist manifolds which are neither complex nor symplectic yet
do admit generalized complex structure.

Since generalized complex manifolds must be almost complex, this
question becomes nontrivial first in dimension 4, where we are
fortunate to have obstructions to the existence of complex and
symplectic structures coming from Seiberg--Witten theory. For
example, a simply-connected complex or symplectic 4-manifold with
$b_+\geq 3$ must have a nonzero Seiberg--Witten
invariant~\cite{Ta94}.

Each tangent space of a generalized complex manifold has a
distinguished subspace equipped with a symplectic form and
transverse complex structure; the transverse complex dimension is
called the \emph{type}, a local invariant of the geometry which may
vary along the manifold.  We show that in 4 dimensions, a connected
and nondegenerate type change locus must be a smooth 2-torus, which
also inherits a complex structure, i.e. it must be a nonsingular
elliptic curve.

We then introduce a surgery for \gcms\ which is a particular case of
the \emph{$C^{\infty}$ logarithmic transformation} introduced by
Gompf and Mrowka \cite{GoMo93}. This surgery modifies a
neighbourhood of a symplectic 2-torus with trivial normal bundle in
a symplectic 4-manifold, producing a new manifold endowed with a
\gcs\ with type change along a 2-torus.  Performing this surgery
along a fiber of an elliptic $K3$ surface, we obtain a \gcs\ on $3
\C P^2 \# 19 \overline{\C P^2}$, a manifold with vanishing
Seiberg--Witten invariants \cite{Wi94,Ta94}.

We thank Tomasz Mrowka for advice which led to our final example. We
also thank Nigel Hitchin for helpful conversations.

\section{Generalized complex structures}\label{sec:gcss}

In this section we recall the definition and basic examples of
generalized complex structures, following~\cite{Gu03}.

Given a closed 3-form $H$ on a manifold $M$, we define the {\it
Courant bracket} of sections of the sum $T\oplus T^*$ of the tangent
and cotangent bundles by
$$[X+\xi,Y+\eta]_H = [X,Y] + \mc{L}_X \eta -\mc{L}_Y\xi -\frac{1}{2}d(\eta(X) - \xi(Y)) + i_Y i_X H.$$
The bundle $T\oplus T^*$ is also endowed with a natural symmetric
pairing of signature $(n,n)$:
$$ \IP{X+\xi,Y+\eta} = \frac{1}{2} (\eta(X) + \xi(Y)).$$

\begin{definition}
A {\it \gcs} on a manifold with closed 3-form  $(M,H)$ is a complex
structure on the bundle $T\oplus T^*$ which preserves the natural
pairing and whose $+i$-eigenspace is closed under the Courant
bracket.
\end{definition}

A \gcs\ can be fully described in terms of its $+i$-eigenspace $L$,
which is a maximal isotropic subspace of $T_{\C} \oplus T^*_{\C}$
satisfying $L \cap \overline{L} = \{0\}$.  Alternatively, it can be
described using differential forms.  Recall that the exterior
algebra $\wedge^{\bullet}T^*$ carries a natural spin representation
for the metric bundle $T \oplus T^*$; the Clifford action of $X+\xi
\in T \oplus T^*$ on $ \rho \in \wedge^{\bullet}T^*$ is
$$ (X+\xi) \cdot \rho  = i_X \rho + \xi \wedge \rho.$$
The subspace $L\subset T_{\C} \oplus T_{\C}^*$ annihilating a spinor
$\rho \in \wedge^{\bullet}T_{\C}^*$ is always isotropic.  If $L$ is
maximal isotropic, then $\rho$ is called a {\it pure spinor} and it
must have the following algebraic form at every point:
\begin{equation}\label{eq:pure}
\rho =  e^{B+ i \go}\wedge \gO,
\end{equation}
where $B$ and  $\go$ are real 2-forms and $\gO$ is a decomposable
complex form.  Pure spinors annihilating the same space must be
equal up to rescaling, hence a maximal isotropic $L \subset T_{\C}
\oplus  T^*_{\C}$ may be uniquely described by a line bundle
$K\subset \wedge^{\bullet}T_{\C}^*$.
\begin{definition}
Given a \gcs\ \J, the line bundle $K
\subset\wedge^{\bullet}T^*_{\C}$ annihilating its $+i$-eigenspace
is the {\it canonical bundle} of $\J$.
\end{definition}
Note that the condition $L \cap \overline{L} = \{0\}$ at $E$ over a point   $p \in M$, is equivalent
to the requirement that
\begin{equation}\label{eq:LcapLbar}
\gO \wedge \overline{\gO} \wedge \go^{n-k} \neq 0
\end{equation}
for any generator $\rho = e^{B+ i \go}\wedge \gO$ of $K$ over $p$, where $k = \deg(\gO)$ and $2n = \dim(M)$. Therefore at
each point of a generalized complex manifold, $\ker
\gO\wedge\overline{\gO}$ is a subspace of the real tangent space
with induced symplectic structure and transverse complex structure.
\begin{definition}
Let $\J$ be a \gcs\ and $e^{B+i \go} \wedge \gO$ a generator of its canonical bundle at a point $p$. The {\it type} of $\J$ at $p$ is the degree of $\gO$. 
\end{definition}
We remark that while the type of a \gcs\ may jump along loci in the
manifold, its parity must remain constant on connected components of
$M$ (see \cite{Gu03}).

Finally, the Courant integrability of $L$ is equivalent to the
requirement that, for any local generator $\rho\in C^\infty(K)$, one
has
\begin{equation}\label{eq:integrability}
d \rho + H\wedge\rho  = v\cdot \rho
\end{equation}
for some section $v \in C^{\infty}(T_{\C} \oplus T^*_{\C})$.  In
summary, a \gcs\ may be specified by a line sub-bundle $K\subset
\wedge^\bullet T^*_\C$ whose local generators satisfy
\eqref{eq:pure}, \eqref{eq:LcapLbar} and \eqref{eq:integrability}.

\begin{ex}{ex:complex} Let $(M^{2n},I)$ be a complex manifold. Then
the following operator on $T\oplus T^*$ is a generalized complex
structure:
$$ \J_I  = \begin{pmatrix} -I & 0 \\ 0 & I^* \end{pmatrix}$$
The $+i$-eigenspace of $\J_I$ is $T^{0,1} \oplus T^{*{1,0}}$, which
annihilates the canonical bundle $K=\wedge^{n,0}T^*$ and is
therefore of type $n$.
\end{ex}

\begin{ex}{ex:symplectic}
Let $(M,\go)$ be a symplectic manifold. Then
$$\J_{\go} = \begin{pmatrix}0& -\go^{-1} \\ \go & 0 \end{pmatrix}$$
is a generalized complex structure with $+i$-eigenspace $\{X -
i\go(X):X \in T_{\C}M\}$ and canonical bundle generated by the
differential form $e^{i \go}$.  Symplectic structures, therefore,
have type zero.
\end{ex}

\begin{ex}{ex:B-field}
A real closed 2-form $B$ gives rise to an orthogonal transformation
of $T\oplus T^*$ via $X+\xi\mapsto X+\xi+i_X B$.  This
transformation, called a \emph{B-field transform}, preserves the
Courant bracket, and hence it acts by conjugation on any given \gcs\
$\J$ on $M$, producing a new one. The induced action on the
canonical bundle is simply $K\mapsto e^B\wedge K$.

If $B$ is not closed, then it induces an isomophism between the
$H$-Courant bracket and the $H+dB$-Courant bracket.  In particular,
if $[H] =0\in H^3(M,\R)$, the bracket $[,]_{H}$ is isomorphic to
$[,]_{0}$ by the action of a nonclosed  2-form.
\end{ex}

In the next example, we demonstrate that the type of a \gcs\ may not
be constant; it jumps from type $0$ to type $2$ along a codimension
2 submanifold.
\begin{ex}{ex:local model}{\em (Local model)}
Consider $\C^2$ with complex coordinates $z_1,z_2$.  The
differential form
$$\rho=  z_{1} + dz_{1}\wedge dz_{2}$$
is equal to $dz_1\wedge dz_2$ along the locus $z_1=0$, while away
from this locus it can be written as
\begin{equation}\label{eq:complex local form}
\rho = {z_1}\exp(\tfrac{dz_1 \wedge dz_2}{z_1}).
\end{equation}
Since it also satisfies $d\rho = -\partial_2\cdot\rho$, we see that
it generates a canonical bundle $K$ for a generalized complex
structure which has type $2$ along $z_1=0$ and type $0$ elsewhere.

Observe that this structure is invariant under translations in the
$z_2$ direction, hence we can take a quotient by the standard $\Z^2$
action to obtain a generalized complex structure on the torus
fibration $D^2 \times T^2$, where $D^2$ is the unit disc in the
$z_1$-plane. Using polar coordinates, $z_1 = r e^{2\pi i\theta_1}$,
the canonical bundle is generated, away from the central fibre, by
\begin{align*}
\exp(B + i \go) &= \exp(d \log r + i d\theta_1)\wedge (d \theta_2 + i d\theta_3) \\
& = \exp(d\log r \wedge d \theta_2 - d\theta_1 \wedge d \theta_3  + i(d\log r \wedge d\theta_3 +d \theta_1 \wedge d \theta_2)),
\end{align*}
where $\theta_2$ and $\theta_3$ are coordinates for the 2-torus with
unit periods. Away from $r=0$, therefore, the structure is a
$B$-field transform of a symplectic structure $\omega$, where
\begin{equation}
\begin{aligned}\label{eq:local model}
B & = d\log r \wedge d \theta_2 - d\theta_1 \wedge d \theta_3\\
\omega & = d\log r \wedge d\theta_3 +d \theta_1 \wedge d \theta_2.
\end{aligned}
\end{equation}
The type jumps from $0$ to $2$ along the central fibre $r=0$,
inducing a complex structure on the restricted tangent bundle, for
which the tangent bundle to the fibre is a complex sub-bundle. Hence
the type change locus inherits the structure of a smooth elliptic
curve with Teichm\"uller parameter $\tau=i$.
\end{ex}

\begin{ex}{logarithmic}
Endow $D^2\times T^2$ with the \gcs\ of Example \ref{ex:local model}
and consider the action of $\Z_m$   given in polar coordinates  by
$$(r, \theta_1,\theta_2,\theta_3) \mapsto (r, \theta_1 + 1/m, \theta_2 + k/m, \theta_3), $$
where $k$ is co-prime with $m$. This action extends to the fiber
over $r=0$, has no fixed points and preserves the \gcs.  Hence the
quotient, which is a singular $T^2$ fibration with multiple central
fibre, has a \gcs. Away from the central fibre, the coordinates $(
r',\theta_1',\theta_2',\theta_3')=(r^m,m\theta_1,
\theta_2-k\theta_1, \theta_3)$ are well-defined, and the generalized
complex structure is generated by $\exp(B + i\omega)$, where
\begin{equation}
\begin{aligned}\label{eq:logarithmic}
B & = d\log r' \wedge (d \theta_2' +\tfrac{k}{m}d\theta_1') - \tfrac{1}{m}d\theta_1' \wedge d \theta_3'\\
\omega & =\tfrac{1}{m}( d\log r' \wedge d\theta_3' +d \theta_1'
\wedge d \theta_2').
\end{aligned}
\end{equation}
Note that the symplectic form is a rescaling of that in
equation~\eqref{eq:local model}.  As in the previous example, the
central fibre obtains a complex structure.
\end{ex}

\section{The type-changing locus}\label{sec:type change locus}

In the last two examples, the type of the generalized complex
structure jumped from 0 to 2 along a 2-torus, which then inherited a
complex structure.  We now show that this happens generically in
four dimensions.

Recall that a \gcm\ has a canonical bundle $K\subset
\wedge^{\bullet}T^*_{\C}$, so the projection from
$\wedge^{\bullet}T^*_{\C}$ onto $\wedge^{0}T^*_{\C} =\C $ determines
a canonical section $s$ of $K^*$.  For a 4-dimensional manifold, the
type of a generalized complex structure jumps from 0 to 2 precisely
when this section vanishes.
\begin{definition}
A point $p$ in the type-changing locus of a \gcs\ on a 4-manifold is
{\it nondegenerate} if it is a nondegenerate zero of $s\in
C^\infty(K^*)$.
\end{definition}
\begin{theo}{t:local form}
The following hold for a 4-dimensional \gcm:
\begin{enumerate}
\item A nondegenerate point in the type-changing locus has a neighbourhood in which the type changes along a smooth 2-manifold with induced complex structure.
\item\label{item:two} A compact connected component of the  type-changing locus whose points are nondegenerate must be a smooth elliptic curve.
\end{enumerate}
\end{theo}
\begin{proof}
To prove the first claim, let $\rho = \rho_0+ \rho_2 +\rho_4$, with
$\deg(\rho_i)=i$, be a nonvanishing local section of $K$ around a
type-changing point $p$. Then $s(\rho) = \rho_0$ and nondegeneracy
implies that $d\rho_0:T_pM \into \C$ is onto. The implicit function
theorem  implies that the zeros of $\rho_0$ near $p$ form a 2
dimensional manifold. According to equations \eqref{eq:pure} and
\eqref{eq:LcapLbar}, $\rho_2$ induces a complex structure on $T_pM$
for which it generates the canonical line $\wedge^{2,0}T_p^*M$. The
integrability condition~\eqref{eq:integrability} states that
$$ d\rho_0 = i_X \rho_2 \in T^{*1,0}_pM,$$
for some $X\in C^{\infty}(T_{\C}M)$.  Therefore $d\rho_0(T^{0,1}_p
M) =0$, showing that the zero set of $\rho_0$ has a complex
structure.

To prove \eqref{item:two}, let $\Sigma \into M$ be a compact
connected component of the type-changing locus  with its induced
complex structure. Then, since  $ds\in C^\infty(T_\C
^*M|_{\Sigma}\tensor K^*|_{\Sigma})$ vanishes on vectors tangent to
$\Sigma$, nondegeneracy implies that $ds$  is a nowhere vanishing
section of $N^*\tensor K|_{\Sigma}^*$, where $N^*$ is the conormal
bundle. In particular, $N^* \cong K|_{\Sigma}$. Since $\J$  is
complex over $\Sigma$, we have an adjunction formula relating the
canonical bundle $K_\Sigma$ of the complex curve $\Sigma$, with the
canonical bundle $K$ restricted to $\Sigma$:
 $$ K|_{\Sigma} \cong K_{\Sigma} \tensor N^*,$$
showing that $K_{\Sigma}$ is trivial and $\Sigma$ is an elliptic
curve, as required.
\end{proof}

\section{The surgery}\label{sec:surgery}

In this section we introduce a surgery for 4-manifolds with \gcs\
which removes a neighborhood of a symplectic 2-torus and replaces it
by a neighborhood of a torus where the \gcs\ changes type, as in
Example \ref{ex:local model}. This surgery  is an example of a {\it
$C^{\infty}$ logarithmic transformation} as defined by Gompf and
Mrowka \cite{GoMo93}, which we now recall.

Let $ T \hookrightarrow M $ be a 2-torus with trivial normal bundle
in a 4-manifold, and let $U\cong D^2\times T^2$ be a tubular
neighborhood. A {\it $C^{\infty}$ logarithmic transform} of $M$ is a
manifold $\tilde{M}$ obtained by removing $U$ and replacing it with
$D^2\times T^2$, glued in by a diffeomorphism $\psi: S^1 \times
T^2\into \del U$:
$$ \tilde{M}  = (M\backslash U) \cup_{\psi}(D^2 \times T^2).$$
The {\it multiplicity} of this transformation is the degree of the
map $\pi \circ \psi:S^1\times {point}\into \del D^2$, where
$\pi:U\into D^2$ is the first projection.

\begin{theo}{t:surgery}
Let $(M,\sigma)$ be a symplectic 4-manifold, $T \hookrightarrow M$
be a symplectic 2-torus with trivial normal bundle and tubular
neighbourhood $U$.  Let $\psi:S^1\times T^2\into \del U\cong
S^1\times T^2$ be the map given on standard coordinates by
$$ \psi(\theta_1,\theta_2,\theta_3) = (\theta_3,\theta_2,- \theta_1).$$
Then the multiplicity zero $C^{\infty}$ logarithmic transform of $M$
along $T$,
$$\tilde{M} = M \backslash U \cup_{\psi}
 D^2 \times T^2,$$ admits a \gcs\ with type change along a 2-torus,
and which is integrable with respect to a 3-form $H$, such that
$[H]$ is the Poincar\'e dual to the circle in $S^1\times T^2$
preserved by $\psi$. If $M$ is simply connected and $[T]\in
H^2(M,\Z)$ is $k$ times a primitive class, then $\pi_1(\tilde{M}) =
\Z_k$.
\end{theo}
\begin{proof}
By Weinstein's neighbourhood theorem~\cite{Wein71}, the neighborhood
$U$ is symplectomorphic to $D^2\times T^2$ with standard sympletic
form:
$$\sigma = \frac{1}{2}d \tilde{r}^2 \wedge d \tilde{\theta}_1 +d\tilde{\theta_2}\wedge
d\tilde{\theta_3}.$$ Now consider the symplectic form $\go$ on
$D^2\backslash\{0\} \times T^2$ from example \ref{ex:local model}:
\[
\omega = d\log r \wedge d\theta_3 +d \theta_1 \wedge d \theta_2.
\]
The map $\psi : (D^2\backslash D^2_{1/\sqrt{e}}\times T^2 , \go)
\into (D^2\backslash \{0\}\times T^2,\sigma)$ given by
$$\psi(r,\theta_1,\theta_2,\theta_3) = (\sqrt{ \log er^2},\theta_3,\theta_2,- \theta_1)$$
is a symplectomorphism.

Let $B$ be the 2-form defined by \eqref{eq:local model} on
$D^2\backslash\{0\} \times T^2$, and choose an extension $\tilde B$
of ${\psi^{-1}}^*B$ to $M\backslash T$.  Therefore $(M\backslash T,
\tilde B + i\sigma)$ is a generalized complex manifold of type 0,
integrable with respect to the $d\tilde B$-Courant bracket.

Now the surgery $\tilde{M} = M\backslash T \cup_{\psi} D^2 \times
T^2$ obtains a generalized complex structure since the gluing map
$\psi$ satisfies $\psi^*(\tilde B +i\sigma) = B + i\omega$, and this
generalized complex structure exhibits type change along the 2-torus
coming from the central fibre of $D^2\times T^2$.  This structure is
integrable with respect to $H=d\tilde B$, which is a globally
defined closed 3-form on $\tilde M$.

The 2-form $\tilde B$ can be chosen so that it vanishes outside a
larger tubular neighbourhood $U'$ of $T$, so that $H=d\tilde B$ has
support in $U'\backslash U$ and has the form
\[
H=f'(r)dr\wedge d\theta_1\wedge d\theta_3,
\]
for a smooth bump function $f$ such that $f|_U=1$ and vanishes
outside $U'$.  Therefore, we see that $H$ represents the Poincar\'e
dual of the circle parametrized by $\theta_2$, as required.

The last claim is a consequence of van Kampen's theorem and that
$H^2(M,\Z)$ is spherical, as $M$ is simply connected.
\end{proof}
\begin{cor}
Since $\tilde B$ can be chosen to have support in a neighbourhood of
the symplectic 2-torus $T$, the surgery above may be performed
simultaneously on a collection of disjoint symplectic 2-tori in $M$.
\end{cor}
Observe that the crucial property of the type-changing \gcs\ on
$D^2\times T^2$  which allows us to perform the surgery is the
behaviour of its symplectic form.  As we saw, this is the same
symplectic form, up to rescaling, as in Example \ref{logarithmic}.
Hence we could, alternatively, use the generalized complex structure
on $(D^2\times T^2)/\Z_m$ described there as a model for the piece
being glued in.

\section{Examples}\label{sec:examples}

\begin{ex}{hopf}
Consider a symplectic 4-manifold $M=\Sigma\times T^2$, where
$\Sigma$ is a symplectic surface and $T^2$ a symplectic 2-torus.
 Performing the surgery from Theorem \ref{t:surgery} along one of the
$T^2$ fibers, we obtain a type-changing \gcs\ on $X^3 \times S^1$,
where $X^3$ is the \emph{twisted connected sum} of the $S^1$-bundles
$\Sigma \times S^1$ and $S^2 \times S^1$, in the language of
\cite{Kov00}.

For example, if $\Sigma = S^2$, we obtain a generalized complex
structure on $S^3 \times S^1$, integrable with respect to a 3-form
$H$ representing a generator for $H^3(S^3,\Z)$. Note that this
manifold does not admit symplectic structure.
\end{ex}

In the final example, we produce a \gcx\ 4-manifold which admits
neither symplectic nor complex structure.  The non-existence of
symplectic structure follows from a result in Seiberg--Witten
theory.

\begin{ex}{SW=0}{(Generalized complex structure on $3 \C P^2 \# 19 \overline{\C P^2}$)}
Consider an elliptically fibred $K3$ surface $M$.  Any smooth
elliptic fibre is a symplectic 2-torus with respect to a K\"ahler
symplectic form, and has trivial normal bundle. Therefore we may
perform our surgery along such a fiber to obtain a \gcm.

The effect of $C^{\infty}$ logarithmic transformations on the $K3$
surface was studied by Gompf and Mrowka in \cite{GoMo93}. Using a
trick of Moishezon \cite{Moi77}, they show that the differentiable
manifold obtained through a transformation of multiplicity zero is
$\tilde{M} = 3 \C P^2 \# 19\overline{\C P^2}$. Since
$H^3(\tilde{M})= \{0\}$, the \gcs\ on $\tilde{M}$ given by Theorem
\ref{t:surgery} has $[H]=0$.

Since $3 \C P^2 \# 19\overline{\C P^2}$ can be expressed as a
connected sum of terms whose intersection forms are not negative
definite, its Seiberg--Witten invariants vanish \cite{Wi94}.
Therefore Taubes' theorem implies it does not admit symplectic
structure \cite{Ta94}.  This, in turn, obstructs the existence of a
complex structure, since Kodaira's theorem~\cite{Kod64} states that
any complex surface with even first Betti number is a deformation of
an algebraic surface.
\end{ex}

\end{document}